\documentclass[12pt]{article} 
\usepackage{latexsym}
\usepackage{amssymb}
\usepackage{euscript}

\def\ie{{\it i.e. }}

\def\={\ = \ }
\def\la{\langle}
\def\ra{\rangle}

\def\be{\setcounter{equation}{\value{theorem}} \begin{equation}}
\def\ee{\end{equation} \addtocounter{theorem}{1}}
\def\bp{{\sc Proof: }}
\def\ep{{}{\hfill $\Box$} \vskip 5pt \par}
\def\bl{\begin{lemma}}
\def\el{\end{lemma}}
\def\bt{\begin{theorem}}
\def\et{\end{theorem}}
\def\bprop{\begin{prop}}
\def\eprop{\end{prop}}

\def\D{{\mathbb D}}

\def\C{{\mathbb C}}

\def\={\ = \ }

\def\l{\lambda}
\def\La{\Lambda}
\def\G{\Gamma}
\def\Ga{\Gamma}
\def\De{\Delta}

\def\hd{H^\infty(\D)}
\def\hb{H^\infty(\D^2)}

\advance\topmargin by -.6in
\advance\textheight by 1.2in
\advance\oddsidemargin by -.5in
\advance\textwidth by 1.05in

\newtheorem{theorem}{Theorem} [section]
\newtheorem{prop}[theorem]{Proposition}
\newtheorem{lemma}[theorem]{Lemma}

\begin{document}
\setlength{\baselineskip}{21pt}
\title{ The three point Pick problem on the Bidisk}
\author{Jim Agler
\thanks{Partially supported by the National Science Foundation}\\
U.C. San Diego\\
La Jolla, California 92093
\and
John E. M\raise.5ex\hbox{c}Carthy\\
Washington University\\
St. Louis, Missouri 63130}
\date{}

\bibliographystyle{plain}

\maketitle


\baselineskip = 18pt

\setcounter{section}{-1}
\section{Introduction}
The original {\it Pick problem} is to determine, given $N$ points
$\l_1, \dots, \l_N$ in the unit disk $\D$ and $N$ complex numbers
$w_1, \dots , w_N$, whether there exist a function $\phi$ in the
closed unit ball of $\hd$ (the space of bounded analytic functions
on $\D$) that maps each point $\l_i$ to the corresponding value
$w_i$. This problem was solved by G.~Pick in 1916 \cite{pi16}, who
showed that a necessary and sufficient condition is that the Pick
matrix
$$
\left( \frac {1 - \bar w_i w_j}{1 - \bar \l_i \l_j}
\right)_{i,j=1}^N
$$
be positive semi-definite.

It is well-known that if the problem is {\it extremal}, \ie the
problem can be solved with a function of norm one but not with a
function of any smaller norm, then the Pick matrix is singular, and
the corresponding solution is a unique Blaschke product, whose
degree equals the rank of the Pick matrix \cite{gar81,foi-fra}.

In \cite{ag1}, the first author extended Pick's theorem to the
space $\hb$, the bounded analytic functions on the 
bidisk; see also \cite{colwer94,baltre98,agmc_bid}.
It was shown in
\cite{agmc_bid} that if the problem has a solution, then it has a
solution that is a rational inner function;
however the qualitative properties of general solutions are not fully
understood. The example $\l_1 = (0,0),\, \l_2 = (\frac 12, \frac
12),\, w_1 = 0,\, w_2 = \frac 12$ shows that even extremal problems
do not always have unique solutions. 

The two point Pick problem on the bidisk is easily analyzed.
It can be solved if and only if the Kobayashi distance between
$\l_1$ and $\l_2$ is greater than or equal to the hyperbolic
distance between $w_1$ and $w_2$. On the bidisk, the Kobayashi
distance is just the maximum of the hyperbolic distance between the
first coordinates, and the hyperbolic distance between the second
coordinates. A pair of points in $\D^2$ is called {\it balanced} if the
hyperbolic distance between their first coordinates equals the
hyperbolic distance between their second coordinates.

The two point Pick problem has a unique solution if and only if
the Kobayashi distance between $\l_1$ and $\l_2$ exactly equals the
hyperbolic distance between $w_1$ and $w_2$, and moreover $(\l_1,\l_2)$
is not balanced. In this case the solution is a M\"obius map in the
coordinate function in which the Kobayashi distance is attained.
If the distance between $\l_1$ and $\l_2$ equals the distance
between $w_1$ and $w_2$, but the pair $(\l_1, \l_2)$ is balanced,
then the function $\phi$ will be uniquely determined on the
geodesic disk passing through $\l_1$ and $\l_2$, but will not be
unique off this disk.
(For the example $\l_1 = (0,0),\, \l_2 = (\frac 12, \frac
12),\, w_1 = 0,\, w_2 = \frac 12$, on the diagonal $\{ (z,z) \}$ we
must have $\phi(z,z) = z$; but off the diagonal any convex
combination of the two coordinate functions $z^1$ and $z^2$ will
work).

It is the purpose of this article to examine the three point Pick
problem on the bidisk.
Our main result is the following:

\bt
The solution to an extremal non-degenerate three point problem on
the bidisk is
unique.
The solution is given by a rational inner function of degree $2$.
There is a formula for the solution in terms of two uniquely
determined rank one matrices.
\et

In the next section, we shall define precisely the terms
``extremal'' and ``non-degenerate'', but roughly it means that the
problem is genuinely two-dimensional, is really a $3$ point problem
not a $2$ point problem, and the minimal norm of a solution is $1$.

\section{Notation and Preliminaries}

We wish to consider the $N$ point Pick interpolation problem
\setcounter{equation}{\value{theorem}}
\begin{eqnarray}
\phi(\l_i) &=& w_i, \qquad i = 1,\dots,N \nonumber \\
\mbox{and} && 
\label{pic}
\\
\| \phi \|_{\hb} &\leq& 1 .
\nonumber
\end{eqnarray}
\addtocounter{theorem}{1}
We shall say that a solution $\phi$ to (\ref{pic})
is an {\it extremal solution} if
$\| \phi \| = 1$, and no solution has a smaller norm.

For a point $\l$ in $\D^2$, we shall use superscripts to denote
coordinates:
$$
\l \= (\l^1,\l^2) .
$$
Let $W, \, \La^1$ and $\La^2$  denote the $N$-by-$N$ matrices
\begin{eqnarray*}
W &\=& \left( 1 - \bar w_i w_j \right)_{i,j=1}^N \\
\La^1 &\=& \left( 1 - \bar \l^1_i \l^1_j \right)_{i,j=1}^N \\
\La^2 &\=& \left( 1 - \bar \l^2_i \l^2_j \right)_{i,j=1}^N .
\end{eqnarray*}

A pair $\Ga,\De$ of $N$-by-$N$ positive semi-definite
matrices is called {\it permissible}
if
\be
\label{rea}
W \= \La^1 \cdot \Ga \ + \ \La^2 \cdot \De .
\ee
Here $\cdot$ denotes the Schur or entrywise product:
$$
(A \cdot B )_{ij} \ := \ A_{ij} B_{ij} .
$$
The main result of \cite{ag1}
is that the problem (\ref{pic}) has a solution if and only if there
is a pair $\Ga,\De$ of permissible matrices.

A {\it kernel} $K$ on $\{ \l_1, \dots , \l_N \} \times 
\{ \l_1, \dots , \l_N \}$ is a positive definite $N$-by-$N$ matrix  
$$
K_{ij} = K(\l_i, \l_j) .
$$
We shall call the kernel $K$ {\it admissible} if 
\begin{eqnarray*}
\La^1 \cdot K \ &\geq& 0 \\
\mbox{and} && \\
\La^2 \cdot K \ &\geq& 0 .
\end{eqnarray*}

If the problem (\ref{pic}) has a solution and $K$ is an admissible
kernel, then (\ref{rea}) implies that $K \cdot W \geq 0$.
We shall call the kernel $K$ {\it active} if it is admissible and
$K \cdot W$ has a non-trivial null-space.
Notice that all extremal problems have an active kernel.

If one can find a pair of permissible matrices one of which is $0$,
then the Pick problem is really a one-dimensional problem because
one can find a solution $\phi$ that depends only on one of the
coordinate functions. If this occurs, we shall call the problem
{\it degenerate}; otherwise we shall call it {\it non-degenerate}.

\section{The three point problem}

We wish to analyze extremal solutions to three point Pick problems.
Fix three points $\l_1, \l_2, \l_3$ in $\D^2$, and three numbers
$w_1,w_2,w_3$. Let notation be as in the previous section.
We shall make the following assumptions throughout
this section:

(a) The function $\phi$ is an extremal solution to the Pick problem
of interpolating $\l_i$ to $w_i$, where $i$ ranges from $1$ to
$3$.

(b) The function $\phi$ is not an extremal solution to any of the
three two point Pick problems mapping two of the $\l_i$'s to the
corresponding $w_i$'s.

(c) The three point problem is non-degenerate.

\bl
If $K$ is admissible, then 
 rank$(K\cdot W) \, > \, 1$.
\label{lem1}
\el

\bp
Suppose $(\Ga,\De)$ is permissible. By (\ref{rea}), we have
$$
K \cdot W \= K \cdot \La^1 \cdot \Ga + K \cdot \La^2 \cdot \De.
$$
If  rank$(K\cdot W) \, = \, 1$, then either $\Ga = 0$
(which violates (c)), or there exists $t>0$ such that
$$
K \cdot \La^1 \cdot \Ga \= t K \cdot W .
$$
But then $(\frac 1t \Ga, 0)$ is permissible, violating assumption
(c).
\ep

\bl
\label{lem2}
If $K$ is an admissible kernel with a non-vanishing column, then 
rank$(K\cdot \La^1) \, \geq \, 2$ and rank$(K\cdot \La^2) \, \geq
\, 2$.
\el
\bp
Suppose that rank$(K\cdot \La^1) \,= \, 1$. 
As no entry of $\La^1$ 
can be $0$, and some column of $K$ is non-vanishing, there is a
column of $K\cdot \La^1$ that is non-vanishing. As 
$K\cdot \La^1$ is self-adjoint and rank one and has non-zero
diagonal entries, the other two columns of $K\cdot \La^1$
must be non-zero multiples of this non-vanishing column.
So $Q := K\cdot \La^1$ is a positive rank one matrix with no zero
entries, and $K$ has no zero entries.

So
\begin{eqnarray*}
\La^2 &\=&  \left( \frac 1K \right) \cdot \left( K \cdot \La^2
\right) \\
&\=& \left( \frac 1Q  \cdot \La^1 \right) \cdot \left( K \cdot \La^2
\right) \\
&\=& \left( \frac 1Q  \cdot K \cdot \La^2 \right) \cdot  \La^1,
\end{eqnarray*}
where by $\frac 1K$ and $\frac 1Q$ is meant the entrywise
reciprocal.
Now $K \cdot \La^2 $ is positive by hypothesis, and 
$\frac 1Q$ is positive because $Q$ is rank one and non-vanishing;
moreover the Schur product of two positive matrices is positive
\cite[Thm 5.2.1]{horjoh91}.
Therefore $(\Ga + \De \cdot  \frac 1Q  \cdot K \cdot \La^2, 0)$
is permissible, which violates assumption (c).
\ep

\bl
If $K$ is an active kernel, it has a non-vanishing column.
\label{lem3}
\el
\bp
By assumption (b), we cannot have both $K(\l_1,\l_2)=0$ and
$K(\l_1,\l_3) = 0$; for then $K$ restricted to $\{\l_2,\l_3\}
\times \{\l_2,\l_3\}$ would be an active kernel for the two point
problem on $\l_2,\l_3$, and so any solution to the two point problem
would have norm at least one, so $\phi$ would be an extremal
solution to the two point problem.

If neither of $K(\l_1,\l_2)$ or $K(\l_1,\l_3)$ are $0$, we are
done. So assume without loss of generality that the first is
non-zero and the second equals zero. But then $K(\l_2,\l_3)$ cannot
equal zero, for then $K$ restricted to $\{\l_1,\l_2\}
\times \{\l_1,\l_2\}$ would be active, violating assumption (b).
Thus we can conclude that the second column of $K$ is
non-vanishing.
\ep

\bl
If $(\Ga,\De)$ is a permissible pair, then rank$(\Ga) \, = \,
1\, =$ rank$(\De)$.
\label{lem4}
\el

\bp
Let $K$ be an active kernel. Then $K \cdot W$ is rank 2, and
annihilates some vector
$$
\vec{\gamma} \= \left(
\begin{array}{c}
\gamma_1 \\
\gamma_2 \\
\gamma_3 
\end{array}
\right) .
$$
Moreover, by assumption (b), none of the entries of $\vec{\gamma}$
are $0$.

Suppose rank$(\Ga) \, > \, 1$. We have
$$
K \cdot W \= K \cdot \La^1 \cdot \Ga + K \cdot \La^2 \cdot \De .
$$
As $ K \cdot \La^1$ has non-zero diagonal terms, Oppenheim's
theorem \cite[Thm 7.8.6]{horjoh85} guarantees that 
rank$(K \cdot \La^1 \cdot \Ga) \, \geq \,$
rank$(\Ga)$.
As $K\cdot W$ has rank 2, and $K \cdot \La^2 \cdot \De \geq 0$, we
must have rank$(\Ga) \, = \, 2$.
Write
$$
\G \= \vec{u} \otimes \vec{u} + \vec{v} \otimes \vec{v},
$$
where $\vec{u}$ and $\vec{v}$ are not collinear;  if
$$
\vec{u} \= \left(
\begin{array}{c}
u_1 \\
u_2 \\
u_3
\end{array}
\right) ,
$$
then $\vec{u} \otimes \vec{u}$ denotes the matrix
$$
\left( \vec{u} \otimes \vec{u} \right)_{ij} \= u_i \bar u_j .
$$

Let 
$$K \cdot \La^1 = \vec{w} \otimes \vec{w} + \vec{x} \otimes
\vec{x}
$$ if $K \cdot \La^1$ is rank two, and
$$
K \cdot \La^1 = \vec{w} \otimes \vec{w} + \vec{x} \otimes
\vec{x} + \vec{y} \otimes \vec{y}
$$ if it is rank three.

Notice that $\left( K \cdot \La^1 \cdot \Ga \right) \vec{\gamma} =
0$, because $K \cdot \La^1 \cdot \Ga$ is positive and
\begin{eqnarray*}
\la \left( K \cdot \La^1 \cdot \Ga \right) \vec{\gamma},
\vec{\gamma} \ra &\=&  
- \la \left( K \cdot \La^2 \cdot \De \right) \vec{\gamma},
\vec{\gamma} \ra \\
&\leq& 0 .
\end{eqnarray*}

Therefore all $4$ of $(\vec{u} \otimes \vec{u}) \cdot
(\vec{w} \otimes \vec{w}),\ (\vec{u} \otimes \vec{u}) \cdot
(\vec{x} \otimes \vec{x}),\
(\vec{v} \otimes \vec{v}) \cdot
(\vec{w} \otimes \vec{w}),\ (\vec{v} \otimes \vec{v}) \cdot
(\vec{x} \otimes \vec{x})$ annihilate $\vec{\gamma}$.
Therefore
\begin{eqnarray*}
\sum_{j=1}^3 \bar u_j \bar w_j \gamma_j &\=& 0\\
&\=& \sum_{j=1}^3 \bar u_j \bar x_j \gamma_j\\
&\=& \sum_{j=1}^3 \bar v_j \bar w_j \gamma_j\\
&\=& \sum_{j=1}^3 \bar v_j \bar x_j \gamma_j
\end{eqnarray*}
Therefore the vectors 
$$
\left(
\begin{array}{c}
\overline{w_1}{ \gamma_1}  \\
\overline{w_2}{ \gamma_2} \\
\overline{w_3}{ \gamma_3}
\end{array}
\right)
\quad
\mbox{and}
\quad
\left(
\begin{array}{c}
\overline{x_1 }{\gamma_1}  \\
\overline{x_2 }{\gamma_2} \\
\overline{x_3 }{\gamma_3}
\end{array}
\right)
$$
are both orthogonal to both $\vec{u}$ and $\vec{v}$, and therefore
are collinear (since $\vec{u}$ and $\vec{v}$ span a two-dimensional
subspace of $\C^3$).
As none of the entries of $\vec{\gamma}$
are $0$, it follows that $\vec{w}$ and $\vec{x}$ are collinear.
Therefore rank$(K\cdot \La^1) \, = \, 1$, contradicting
Lemmata~(\ref{lem2}) and (\ref{lem3}).
\ep

\bl
The matrices $\Ga$ and $\De$ are unique.
\label{lem5}
\el

\bp
If both $(\Ga_1,\De_1)$ and $(\Ga_2,\De_2)$
were permissible, then $\displaystyle ( \frac 12 (\Ga_1 + \Ga_2),
\frac 12 (\De_1 + \De_2) )$ would also be permissible.
As all permissible matrices are rank one by Lemma~\ref{lem4},
it follows that $\Ga_1$ and $\Ga_2$ are constant multiples of each
other, and so are $\De_1$ and $\De_2$.

So suppose
\begin{eqnarray*}
W &\=& \La^1 \cdot \Ga + \La^2 \cdot \De \\
\mbox{and}&&\\
W &\=& \La^1 \cdot t_1 \Ga + \La^2 \cdot t_2 \De ,
\end{eqnarray*}
where both $t_1, t_2$ are positive, one is less than $1$, and 
the other is bigger than $1$.
Then
$$
(1-t_1) \La^1 \cdot \Gamma + 
(1-t_2) \La^2 \cdot \De \= 0.
$$
Assume without loss of generality that 
$t_1 < 1 < t_2$.
Then $\displaystyle (\frac{t_2 - t_1}{1-t_1} \Ga , 0)$
is permissible, which contradicts Assumption (c).
\ep

\bt
The solution to an extremal non-degenerate three point problem satisfying
Assumptions (a)-(c) is unique.
It is given by a rational inner function of degree $2$, 
and there is
a formula in terms of $\Ga$ and $\De$.
\et
\bp
We have
\be
\label{eq6}
W \= \La^1 \cdot \Ga + \La^2 \cdot \De .
\ee
Choose vectors $\vec{a}$ and $\vec{b}$ so that $\Ga = \vec{a}
\otimes \vec{a}$ and $\De = \vec{b}
\otimes \vec{b}$.

Choose some point $\l_4$ in $\D^2$, distinct from the first three
points. Let $w_4$ be the value attained at $\l_4$ by some solution
$\phi$ of the three point problem~(\ref{pic}).
Then the four point problem, interpolating $\l_i$ to $w_i$ for $i =
1, \dots, 4$ has a solution, so we can find a pair of $4$-by-$4$
permissible matrices $\tilde \Ga$ and $\tilde \De$ satisfying~(\ref{rea}).
As the restriction of these matrices to 
the first three points 
satisfy (\ref{eq6}), and $\Ga$ and $\De$ are unique by
Lemma~\ref{lem5}, we get that $\tilde \Ga$ and $\tilde \De$
are extensions of $\Ga$ and $\De$.
Therefore we have
\setcounter{equation}{\value{theorem}}
\begin{eqnarray}
\left(
\begin{array}{cccc|c}
&&&&1- \bar w_1 w_4 \\
&W&&& 1 - \bar w_2 w_4\\
&&&&1 - \bar w_3 w_4 \\ \cline{1-4}
*&*&*&\multicolumn{2}{r}{ 1 - |w_4|^2}
\end{array}
\right)
\ &=& \
\left(
\begin{array}{cccc|c}
&&&&g_1\\
&\Ga&&& g_2\\
&&&&g_3 \\ \cline{1-4}
\bar g_1&\bar g_2&\bar g_3&\multicolumn{2}{r}{ g_4}
\end{array}
\right)
\,
\cdot 
\,
\left(
\begin{array}{cccc|c}
&&&&1 - \bar \l_1^1 \l_4^1 \\
&\La^1&&& 1 - \bar \l_2^1 \l_4^1\\
&&&&1 - \bar \l_3^1 \l_4^1 \\ \cline{1-4}
*&*&*&\multicolumn{2}{r}{ 1 - |\l_4^1|^2}
\end{array}
\right)
\nonumber
\\
\label{bigp}
&&\\
&&+ \
\left(
\begin{array}{cccc|c}
&&&&d_1\\
&\De&&& d_2\\
&&&&d_3 \\ \cline{1-4}
{\rule[0mm]{0mm}{4mm}
\bar d_1}&\bar d_2&\bar d_3&\multicolumn{2}{r}{ d_4}
\end{array}
\right)
\,
\cdot \,
\left(
\begin{array}{cccc|c}
&&&&1 - \bar \l_1^2 \l_4^2 \\
&\La^2&&& 1 - \bar \l_2^2 \l_4^2\\
&&&&1 - \bar \l_3^2 \l_4^2 \\ \cline{1-4}
*&*&*&\multicolumn{2}{r}{ 1 - |\l_4^2|^2}
\end{array}
\right)
\nonumber
\end{eqnarray}
\addtocounter{theorem}{1}

As $\tilde \Ga$ is positive, it must be that $\displaystyle
\left(
\begin{array}{c}
g_1 \\
g_2 \\
g_3
\end{array}
\right)
$ is in the range of $\Ga$, so 
$$
\left(
\begin{array}{c}
g_1 \\
g_2 \\
g_3
\end{array}
\right) \=
s \vec{a} \= 
s
\left(
\begin{array}{c}
a_1 \\
a_2 \\
a_3
\end{array}
\right)
$$
for some constant $s$.
Similarly,
$$
\left(
\begin{array}{c}
d_1 \\
d_2 \\
d_3
\end{array}
\right) \=
t \vec{b}
\= t
\left(
\begin{array}{c}
b_1 \\
b_2 \\
b_3
\end{array}
\right)
$$
for some $t$.

Let
\begin{eqnarray*}
\vec{v_1} &\=& 
\left(
\begin{array}{c}
(1 - \bar \l_1^1 \l_4^1)a_1 \\
(1 - \bar \l_2^1 \l_4^1)a_2 \\
(1 - \bar \l_3^1 \l_4^1)a_3
\end{array}
\right)
\\
\vec{v_2} &\=&
\left(
\begin{array}{c}
(1 - \bar \l_1^2 \l_4^2)b_1 \\
(1 - \bar \l_2^2 \l_4^2)b_2 \\
(1 - \bar \l_3^2 \l_4^2)b_3
\end{array}
\right)
\\
\vec{v_3} &\=&
\left(
\begin{array}{c}
\bar w_1 \\
\bar w_2 \\
\bar w_3
\end{array}
\right).
\end{eqnarray*}

Looking at the first three entries of the last column of
Equation~(\ref{bigp}), we get
\be
\label{59}
\left(
\begin{array}{c}
1 \\
1 \\
1
\end{array}
\right)
\=
s \vec{v_1} \, + \, t \vec{v_2} \, + \,
w_4 \vec{v_3}.
\ee
Equation~(\ref{59}) has a unique solution for $s,t$ and $w_4$
unless
\be
\label{60}
det 
\left(
\begin{array}{ccc}
(1 - \bar \l_1^1 \l_4^1)a_1 & (1 - \bar \l_1^2 \l_4^2)b_1 & \bar
w_1 \\
(1 - \bar \l_2^1 \l_4^1)a_2 & (1 - \bar \l_2^2 \l_4^2)b_2&
\bar w_2 \\
(1 - \bar \l_3^1 \l_4^1)a_3 & (1 - \bar \l_3^2 \l_4^2)b_3&
\bar w_3
\end{array}
\right)
\= 0 .
\ee

Notice that the determinant in (\ref{60}) is analytic in $\l_4$.
So if there is a single point $\l_4$ for which the determinant
does not vanish, there is an open neighborhood of this point for
which the determinant doesn't vanish. Consequently $w_4$ (and hence
$\phi$) would be determined uniquely on this open set, and hence on
all of $\D^2$.

Suppose the determinant in (\ref{60}) vanished identically.
Then there is a set of uniqueness of $\l_4$'s on which 
Equation~(\ref{59}) can be solved with either $s$ or $t$ equal to
$0$. (If both $s$ and $t$ were uniquely determined, then 
$\vec{v_3}$ would be $0$, violating Assumption (a)).
Without loss of generality, take $t=0$.
Moreover, we can also assume without loss of generality that 
$w_1$ and $w_2$ do not both vanish.

Then one can use the first component of Equation~(\ref{59}) to
solve for $s$, and the second one to get
$$
w_4 \= \frac{(1 - \bar \l_1^1 \l_4^1)a_1 + (1 - \bar \l_2^1
\l_4^1)a_2}
{(1 - \bar \l_1^1 \l_4^1)a_1 \bar w_2 + 
(1 - \bar \l_2^1 \l_4^1)a_2 \bar w_1}
$$
Then $w_4$ is given uniquely as a rational function of degree $1$
of $\l_4^1$, violating both Assumptions (b) and (c).

Therefore we can assume that there is an open set on which
Equation~(\ref{59}) has a unique solution, so by Cramer's rule we
get 
\be
\label{61}
\phi(\l_4) \= w_4 \=
\frac{ det 
\left(
\begin{array}{ccc}
(1 - \bar \l_1^1 \l_4^1)a_1 & (1 - \bar \l_1^2 \l_4^2)b_1 & 1 \\
(1 - \bar \l_2^1 \l_4^1)a_2 & (1 - \bar \l_2^2 \l_4^2)b_2&
1 \\
(1 - \bar \l_3^1 \l_4^1)a_3 & (1 - \bar \l_3^2 \l_4^2)b_3&
1
\end{array}
\right)}
{ det
\left(
\begin{array}{ccc}
(1 - \bar \l_1^1 \l_4^1)a_1 & (1 - \bar \l_1^2 \l_4^2)b_1 & \bar
w_1 \\
(1 - \bar \l_2^1 \l_4^1)a_2 & (1 - \bar \l_2^2 \l_4^2)b_2&
\bar w_2 \\
(1 - \bar \l_3^1 \l_4^1)a_3 & (1 - \bar \l_3^2 \l_4^2)b_3&
\bar w_3
\end{array}
\right)}
.
\ee
Equation~(\ref{61})
gives a formula for $\phi$ that shows that $\phi$ is a rational
function of degree at most $2$, whose second order terms only
involve the mixed product $\l_4^1 \l_4^2$.

To show that $\phi$ is inner, we follow \cite{agmc_bid}.
We can rewrite~(\ref{rea}) as
\be
\label{63}
1 + \bar \l_i^1 \l_j^1 a_i \bar a_j + 
\bar \l_i^2 \l_j^2 b_i \bar b_j
\=
\bar w_i w_j + a_i \bar a_j + b_i \bar b_j.
\ee
Realizing both sides of (\ref{63}) as Grammians, we get that there
exists a $3$-by-$3$ unitary $U$ such that, for $j = 1,2,3$,
\be
\label{64}
U  \left(
\begin{array}{c}
1 \\
\l_j^1 \bar a_j \\
\l_j^2 \bar b_j
\end{array}
\right)
\=
 \left(
\begin{array}{c}
w_j \\
\bar a_j \\
\bar b_j
\end{array}
\right).
\ee
Writing
$$
U \=
\begin{array}{cc}
&
\!\!\!\!\!\!
\C \: \  \C^2 \\
\begin{array}{c}
\C \\
\C^2
\end{array} &
\!\!\!\!\!\!
\left(
\begin{array}{cc}
A & B \\
C & D
\end{array}
\right) ,
\\ &
\end{array}
$$
and letting 
$$
E_\l \= \left(
\begin{array}{cc}
\l^1 & 0 \\
0 & \l^2
\end{array}
\right) ,
$$
we can solve (\ref{64}) to get
$$
w_j \= A + B E_{\l_j} ( 1 - D E_{\l_j})^{-1} C .
$$
So the function 
$$
\psi (\l) \= A + B E_{\l} ( 1 - D E_{\l})^{-1} C
$$
interpolates the original data. Moreover $\psi$ is inner,
 because a calculation
shows that
$$
1- \overline{\psi(\lambda)}{\psi(\lambda)} \=
((1 - D E_\lambda)^{-1} C)^\ast (1 - E_\lambda^\ast E_\lambda) (
(1 - D E_\lambda)^{-1} C),
$$
so $|\psi|$ is less than $1$ on $\D^2$ and equals $1$ on the
distinguished boundary.
By uniqueness, we must have $\psi = \phi$, and hence $\phi$ is
inner.

Finally, we must show that the degree of $\phi$ is exactly two.
This is because an easy calculation shows that a rational function
of degree one
$$
\frac{ c_1 + c_2 z^1 + c_3 z^2}{c_4 + c_5 z^1 + c_6 z^2}
$$
is inner only if it is a function of either just $z^1$ or just
$z^2$,
\ie either both $c_2$ and $c_5$ or both $c_3$ and $c_6$ can be
chosen to be zero. This would violate Assumption (c).
\ep

\section{Finding $\Gamma$ and $\Delta$}

Formula~(\ref{61}) works fine, provided one knows $\Gamma$ and
$\Delta$ (or, equivalently, $a_1, a_2, a_3$ and $b_1,b_2, b_3$).
Lemma~\ref{lem5} assures us that $\Gamma$ and $\Delta$ are unique;
how does one find them?

First, let us make a simplifying normalization. One can pre-compose
$\phi$ with an automorphism of $\D^2$, and post-compose it with an
automorphism of $\D$; so one can assume that $\l_1 = (0,0)$ and
$w_1 = 0$. Write $\l_2 = (\alpha_2, \beta_2)$ and 
$\l_3 = (\alpha_3,\beta_3)$.
Moreover, as $\G_{ij} = a_i \bar a_j $ and $\Delta_{ij} = b_i \bar
b_j$, we can choose $a_1 \geq 0$ and $b_1 \geq 0$; again without
loss of generality we can assume that $b_1 > 0$. Thus we have
\setcounter{equation}{\value{theorem}}
\begin{eqnarray}
\left(
\begin{array}{ccc}
1&1&1 \\
1&1-|w_2|^2& 1 - \bar w_2 w_3\\
1&1 - w_2 \bar w_3  & 1 - |w_3|^2 
\end{array}
\right)
\ &=& \
\left(
\begin{array}{ccc}
a_1^2&a_1 \bar a_2 &a_1 \bar a_3  \\
a_1 a_2 &|a_2|^2&  a_2 \bar a_3 \\
a_1 a_3 & \bar a_2 a_3 & |a_3|^2
\end{array}
\right)
\ \cdot \
\left(
\begin{array}{ccc}
1&1&1 \\
1&1-|\alpha_2|^2& 1 - \bar \alpha_2 \alpha_3\\
1&1 - \alpha_2 \bar \alpha_3 & 1 - |\alpha_3|^2
\end{array}
\right)
\nonumber
\\
\label{jan1}&& \\
&&+
\left(
\begin{array}{ccc}
b_1^2&b_1 \bar b_2 &b_1 \bar b_3  \\
b_1 b_2 &|b_2|^2&  b_2 \bar b_3 \\
b_1 b_3 & \bar b_2 b_3 & |b_3|^2
\end{array}
\right)
\ \cdot \
\left(
\begin{array}{ccc}
1&1&1 \\
1&1-|\beta_2|^2& 1 - \bar \beta_2 \beta_3\\
1&1 - \beta_2 \bar \beta_3 & 1 - |\beta_3|^2
\end{array}
\right)
\nonumber
\end{eqnarray}
\addtocounter{theorem}{1}
Looking at the first column of (\ref{jan1})
we get
\begin{eqnarray*}
b_1 &\=& \frac{1-a_1^2}{\sqrt{1- a_1^2}}\\
b_2 &\=& \frac{1-a_1 a_2 }{\sqrt{1- a_1^2}}\\
b_3 &\=& \frac{1-a_1 a_3 }{\sqrt{1- a_1^2}}
\end{eqnarray*}
Thus we have three equations that {\it uniquely} determine $a_1,
a_2, a_3$:
\begin{eqnarray}
\label{jan11}
(1- a_1^2) (1 - |w_2|^2) &=& (1- a_1^2) |a_2|^2 (1 - |\alpha_2|^2)
+ |1 - a_1 a_2|^2 (1 - |\beta_2|^2)
\\
\label{jan12}
(1- a_1^2)(1- \bar w_2 w_3) &=& (1- a_1^2) a_2 \bar a_3 (1 -
\bar \alpha_2 \alpha_3 ) + (1 - a_1 a_2) (1 - a_1 \bar a_3) (1 -
\bar \beta_2 \beta_3) \\
\label{jan13}
(1- a_1^2)(1- |w_3|^2) &=& (1- a_1^2) |a_3|^2 (1 - |\alpha_3|^2)
+ |1 - a_1 a_3|^2 (1 - |\beta_3|^2)
\end{eqnarray}
\addtocounter{theorem}{3}
Equation~(\ref{jan12}) can be used to solve for $a_3$ as a rational
function of $a_1$ and $\bar a_2$; then one is left with two real
algebraic equations in three real variables, $a_1, \Re(a_2)$ and
$\Im(a_3)$.
Provided the original data is really extremal, this system of
two equations will have a unique solution with $a_1 \geq 0$.
If the original data is not extremal, multiply $w_2$ and $w_3$ by a
positive real number $t$, and choose the largest $t$ for which
equations~(\ref{jan11})--(\ref{jan13}) can be solved. This will
produce an inner function $\phi$ via (\ref{61}); then the function
$\frac 1t \phi$ will be the unique function of minimal norm solving
the original problem.

\addtocounter{theorem}{1}
{\bf Example {\arabic{section}.\arabic{theorem}}}
Let us consider a
very symmetric special case.
Let $\l_1 = (0,0), \, \l_2 = (r,0), \, \l_3 = (0,r),
\, w_1 = 0, \, w_2 = t$ and $w_3 = t$, where $t$ is to be chosen as
large as possible and $r$ is a fixed {\it positive} number.
Then by symmetry, we can assume that $a_1 = b_1 =
\frac{1}{\sqrt{2}}$ and $a_2 = b_3 =
\bar a_2$.

Equations~(\ref{jan11})--(\ref{jan13}) then reduce to:
\begin{eqnarray*}
\frac 12 (1 - t^2) &\=& \frac 12 a_2^2 (1 - r^2) + (1 -
\frac{1}{\sqrt{2}} a_2 )^2 \\
\frac 12 (1 - t^2) &\=&  a_2 (\sqrt{2} - a_2) 
\end{eqnarray*}
Solving, one gets two solutions.
One solution is 
$$
t = \frac{r}{2-r}, \quad a_2 = \frac{\sqrt{2}}{2-r} ; 
$$
the other is
$$
t = \frac{r}{2+r}, \quad a_2 = \frac{\sqrt{2}}{2+r}.
$$
The first of these is clearly the extremal solution, and
formula~(\ref{61}) then gives
$$
\phi(z) \= \frac{z^1 + z^2 - 2 z^1z^2}{2 - z^1 - z^2}
$$
as the extremal solution.

The second solution also corresponds to a pair of rank one matrices
$\Ga$ and $\De$ that satisfy~(\ref{eq6}), even though the problem
is non-extremal. If one plugs in to (\ref{61}) one gets the inner
function
$$
\phi_2(z) \= \frac{z^1 + z^2 + 2 z^1z^2}{2 + z^1 + z^2}.
$$


\end{document}